\documentclass[12pt]{amsart}

\usepackage{amssymb}
\usepackage{amsmath}

\newcommand{\Pf}{{\em Proof}. }
\newcommand{\EPf}{\hfill$\square$}
\newcommand{\Z}{\mbox{$\mathbf Z$}}
\newcommand{\R}{\mbox{$\mathbf R$}}
\newcommand{\C}{\mbox{$\mathbf C$}}
\newcommand{\Q}{\mbox{$\mathbf H$}}

\newcommand{\SU}[1]{\mbox{$\mathbf{SU}(#1)$}}
\newcommand{\U}[1]{\mbox{$\mathbf{U}(#1)$}}
\newcommand{\SP}[1]{\mbox{$\mathbf{Sp}(#1)$}}
\newcommand{\SO}[1]{\mbox{$\mathbf{SO}(#1)$}}

\newcommand{\Spin}[1]{\mbox{$\mathbf{Spin}(#1)$}}
\newcommand{\F}{\mbox{$\mathbf{F}_4$}}
\newcommand{\G}{\mbox{$\mathbf{G}_2$}}

\newtheorem{thm}{Theorem}

\newtheorem{prop}[thm]{Proposition}
\newtheorem{lem}[thm]{Lemma}

\begin{document}

\title{Polar actions on compact rank one symmetric spaces
are taut}
\author{Leonardo Biliotti and Claudio Gorodski}
\address{Dipartimento di Matematica, Universit\`a Politecnica 
delle Marche, Via Brecce Bianche, 60131, Ancona, Italy}
\email{biliotti@dipmat.univpm.it}
\address{Instituto de Matem\' atica e Estat\'\i stica,
         Universidade de S\~ ao Paulo,
         Rua do Mat\~ao, 1010,
         S\~ ao Paulo, SP 05508-090,
         Brazil}
\email{gorodski@ime.usp.br}
\thanks{The second author was supported in part by 
FAPESP and CNPq.}

\date{\today}

\begin{abstract}
We prove that the orbits of a polar action of a compact Lie group 
on a compact rank one symmetric space 
are tautly embedded with respect to $\Z_2$-coefficients.
\end{abstract}

\maketitle

\section{Introduction}

The main result of this paper is the following theorem.

\begin{thm}\label{main}
A polar action of a compact Lie group 
on a compact rank one symmetric space 
is taut.
\end{thm}

A proper isometric action of a Lie group $G$ on a
complete Riemaniann manifold $X$ is called 
\emph{polar} if there exists a connected,
complete submanifold $\Sigma$ that meets
all orbits of~$G$ in such a way that the intersections between $\Sigma$ and
the orbits of $G$ are all orthogonal. Such a submanifold
is called a \emph{section}. It is easily seen that a section
$\Sigma$ is totally geodesic in
$X$. An action admitting a section that is flat in the induced metric is
called \emph{hyperpolar}. See the monographs~\cite{P-T,BCO},
and the survey~\cite{Th4} for recent results and a list
of references. 

A properly embedded submanifold $M$ of a complete
Riemannian manifold $X$ is called \emph{$F$-taut},
where $F$ is a coefficient field, 
if the energy functional $E_q:\mathcal P(X,M\times q)\to\mathbf R$
is a $F$-perfect Morse function for every $q\in X$ that is not a focal point
of $M$, where $\mathcal P(X,M\times q)$ denotes the space of $H^1$-paths
$\gamma:[0,1]\to X$ such that $\gamma(0)\in M$ and 
$\gamma(1)=q$~\cite{GH,TTh2}. 
Recall that $E_q$ is a Morse function
if and only if $q$ is not a focal point of $M$, and 
$\gamma$ is a critical point of $E_q$ if and
only if $\gamma$ is a geodesic perpendicular to $M$ 
at $\gamma(0)$ that is parametrized proportionally to arc length,
in which case the index of $\gamma$ as a critical 
point of $E_q$ is 
equal to the sum of the multiplicities of the
focal points of $M$ along $\gamma$
by the Morse index theorem. 
Write $\mathcal M = \mathcal P(X,M\times q)$
and let $\mathcal M^c$ denote the sublevel set $\{E_q\leq c\}$. 
Since $E_q$ is bounded below and satisfies 
the Palais-Smale condition~\cite{P-S}
(see e.g.~Theorem~A in~\cite{GH}),
the number $\mu_k(E_q|\mathcal M^c)$ of critical points of index 
$k$ of a Morse function $E_q$ restricted to $\mathcal M^c$ 
is finite and the weak Morse
inequalities say that 
$\mu_k(E_q|\mathcal M^c)\geq\beta_k(\mathcal M^c,F)$ for all
$k$ and all $c$ and all $F$, 
where $\beta_k(\mathcal M^c,F)$ denotes the $k$-th
$F$-Betti number of $\mathcal M^c$;
$E_q$ is called \emph{$F$-perfect} if 
$\mu_k(E_q|\mathcal M^c)=\beta_k(\mathcal M^c,F)$ for all $k$ and all 
$c$. 

A proper isometric action of a Lie group
$G$ on a complete Riemannian manifold $X$ is 
called \emph{taut} if all the orbits of $G$ in $X$ are 
taut submanifolds. Herein we will always consider 
$F=\Z_2$ and usually drop the reference to the field. 

A proper isometric action of a Lie group
is taut (resp.~polar) if and only if
the same is true for the restriction of that action to the connected
component of the group. This is clear  
in the case of taut actions and is easy in the
case of polar actions, see Proposition~2.4 in~\cite{HPTT2}.
We can thus assume that we have a connected group
in the statement of Theorem~\ref{main}. 

It follows from
results of Bott-Samelson and Conlon~\cite{B-S,Co}
that a hyperpolar action of a compact Lie group
on a complete Riemannian manifold is taut.
(Hyperpolar actions on 
compact irreducible symmetric spaces
were classified by Kollross~\cite{Kol1}.)   
However, there are many examples of 
polar actions on compact rank one symmetric spaces
that are not hyperpolar 
as shown by the classification of polar
actions on these spaces due 
to Podest\`a and Thorbergsson~\cite{PTh1}. 
On the other hand, it has been conjectured
for some time that the classes of polar 
and hyperpolar actions on 
compact irreducible symmetric spaces 
of higher rank coincide. This has been 
verified for actions on these spaces with a fixed 
point~\cite{Br}, real Grassmannians of rank two~\cite{PTh2},
complex Grassmannians~\cite{BG}, and, more recently, 
this result has also been claimed for  
compact irreducible Hermitian symmetric spaces~\cite{Bil},
and compact irreducible type~I symmetric spaces~\cite{Kol2}. 

The compact rank one symmetric spaces are 
the spheres and the projective spaces over the four normed 
division algebras. Regarding the proof of Theorem~\ref{main},
the cases of spheres 
and real projective spaces are very simple.
In the cases of complex and quaternionic projective 
spaces, we use a lifting argument 
based on results of~\cite{GH,TTh2,HLO},
see Proposition~\ref{prop:lifting}.  
It is only in the case of the Cayley projective plane
$\mathbf{Ca}P^2$ that we resort to the classification
in~\cite{PTh1}. 
According to that result, a polar
action on $\mathbf{Ca}P^2$ has either cohomogeneity one 
so that it is hyperpolar and then it is taut
by the results of Bott-Samelson and Conlon referred to above,
or it has cohomogeneity 
two and is given by one of the following subgroups
of $\F$, the isometry group of $\mathbf{Ca}P^2$:
\begin{equation}\label{gps}
 \Spin8,\quad S^1\cdot\Spin7,\quad\SU2\cdot\SU4,\quad\SU3\cdot\SU3.
\end{equation}
We prove that each one of these actions is taut
by generalizing an argument based on the 
reduction principle that was first used 
in~\cite{G-Th1} to prove that certain 
finite dimensional orthogonal representations are taut,
see Proposition~\ref{red}. 

As is easy to see and we will explain in a moment, 
the classes of taut actions on $S^n$, $\R P^n$ and $\R^{n+1}$
all coincide. (A linear action on an Euclidean space
is usually called a representation.)
On the other hand, Proposition~\ref{prop:lifting}
implies that the class of taut actions on $\C P^n$ (resp.~$\Q P^n$)
corresponds to a certain subclass of taut representations. 
In particular, it is worth noting that 
there are taut actions on $\C P^n$ and $\Q P^n$ which are not polar. 
In fact, Proposition~\ref{prop:lifting} implies that the
representations of $\SO2\times\Spin9$ (resp.~$\U2\times\SP n$,
$\SP1\times\SP n$, where $n\geq2$)
that were shown to be taut in~\cite{G-Th1,G-Th2}
induce taut actions on $\C P^{15}$ (resp.~$\C P^{4n-1}$, $\Q P^{2n-1}$),
and these are not polar by the classification in~\cite{PTh1}.

Part of this work has been completed while
the authors were visiting University of Florence,
for which they wish to thank Fabio Podest\`a 
for his hospitality. 

\section{A lifting argument}\label{lifting}

We begin now the proof of Theorem~\ref{main}.
Consider the unit sphere $S^n$ in Euclidean space 
$\R^{n+1}$.  
A polar action action of $G$ on $S^n$ is the restriction 
of a polar representation of $G$ on $\R^{n+1}$.
On Euclidean space, the classes of polar and 
hyperpolar representations coincide since the sections 
are linear subspaces. It follows from the results of
Bott-Samelson and Conlon quoted in the introduction
that a polar representation is taut.
One easily sees that a compact submanifold of $S^n$ is taut in 
$S^n$ if and only if it is taut in $\R^{n+1}$
(this can be seen by recalling 
that in the cases in which $X=S^n$ or $X=\R^{n+1}$, 
a properly embedded submanifold $M$ of $X$
is taut if and only the squared distance function 
$L_q:M\to\R$, $L_q(x)=d(x,q)^2$ is perfect for generic $q\in X$,
see section~2 in~\cite{TTh2}).
It follows that the action of $G$ on $S^n$ is taut. 

Next we consider the cases of polar actions 
on real, complex or quaternionic
projective spaces. By making use of the 
principal bundles 
\[ \Z_2\to S^n\to\R P^n,\qquad S^1\to S^{2n+1} \to \C P^n,\qquad S^3\to S^{4n+3} \to \Q P^n, \]
we reduce the problem to the already solved
case of spheres via the following lifting argument (see 
Lemma~6.1 in~\cite{HLO}; see also Proposition~4.1
in~\cite{GH} and Theorem~6.12 in~\cite{TTh2};
the case of $\R P^n$ is of course simpler and can be dealt with 
directly, compare~Proposition~6.3, p.~45, in~\cite{Be}).

\begin{lem}(Heintze, Liu, Olmos~\cite{HLO})\label{hlo}
Let $\hat X$ and $X$ be Riemannian manifolds,
$\pi:\hat X\to X$ a Riemannian submersion,
$M$ a properly embedded submanifold of $X$, and $\hat M=\pi^{-1}(M)$.
Then a normal vector $\hat\xi\in\nu\hat X$ is a 
multiplicity-$m$ focal direction of $\hat M$
if and only if $\pi_*\hat\xi$ is a 
multiplicity-$m$ focal direction of $M$.
\end{lem}

\begin{prop}\label{prop:lifting}
Let $\hat X$ and $X$ be Riemannian manifolds, and suppose
there is a proper isometric action of a Lie group $G$ on $\hat X$
such that $\pi:\hat X\to X$ is a principal $G$-fiber bundle
and a Riemannian submersion. 
Let $M$ be a properly embedded submanifold of $X$ and $\hat M=\pi^{-1}(M)$.
Then $\hat M$ is taut in $\hat X$ if and only if $M$ is taut in $X$.
\end{prop}

\Pf Let $\hat q\in\hat X$ and $q=\pi(\hat q)$.
It follows from Lemma~\ref{hlo}
that $\hat q$ is not a focal point of $\hat M$
if and only if $q$ is not a focal point of $M$. 
Consider the energy functionals 
$E_{\hat q}:\mathcal P(\hat X,\hat M\times\hat q)\to\R$
and $E_q:\mathcal P(X,M\times q)\to\R$. 
A critical point of $E_{\hat q}$ is a geodesic
joining $q$ to a point in $\hat M$ that is perpendicular
to $\hat M$, and similarly for the critical points 
of $E_q$. Since a geodesic in $\hat X$ is horizontal
with respect to $\pi$ if and only if it is perpendicular
to one fiber of $\pi$, it follows that the critical points
of $E_{\hat q}$ are horizontal geodesics 
and there is a bijective correspondence between 
critical points of $E_{\hat q}$ and critical points
of $E_q$ given by the projection. By the Morse
index theorem, the index of a geodesic as a critical point
of the energy functional is the sum of the
multiplicities of the focal points along the geodesic.
It follows from this and Lemma~\ref{hlo}
that the projection maps critical points of $E_{\hat q}$ to
critical points of $E_q$ with
the same value of the energy and the same indices. In order to finish
the proof of the lemma, we only need to show that 
$\mathcal P(\hat X,\hat M\times\hat q)^c$
and $\mathcal P(X,M\times q)^c$ have the same 
homotopy type for all $c$. In fact, the map
\[ \Phi:\mathcal P(X,M\times q)\times\mathcal P(G,G\times 1)\to
\mathcal P(\hat X,\hat M\times\hat q), \]
given by $\Phi(x,g)(t)=\hat x(t)g(t)$, where 
$\hat x(t)$ is the horizontal lift 
of $x(t)$ with $\hat x(1)=\hat q$, is a 
diffeomorphism. Since $\mathcal P(G,G\times 1)$ is
contractible, the result follows. \EPf 

\medskip

It follows from Proposition~\ref{prop:lifting} 
that polar actions on $\R P^n$,
$\C P^n$ and $\Q P^n$ are taut. 
In fact, the case of $\R P^n$ is clear, 
and if $G$ acts polarly on $\C P^n$ (resp.~$\Q P^n$),
then $T^1\times G$ (resp.~$\SP1\times G$)
acts polarly (and hence tautly)
on $S^{2n+1}\subset\C^{n+1}$
(resp.~$S^{4n+3}\subset\Q^{n+1}$)
by~\cite{PTh1}. Hence $G$ acts tautly on 
$\C P^n$ (resp.~$\Q P^n$) by Proposition~\ref{prop:lifting}.

It remains only to consider the case of the Cayley
projective plane $\mathbf{Ca}P^2$. In this case, there
is no convenient  Riemannian submersion from a sphere,
so we need to use a different argument. 
According to~\cite{PTh1} and the discussion in the introduction, 
we need to prove that the four polar actions 
on $\mathbf {Ca}P^2$ given by
the groups~$G$ listed in~(\ref{gps}) are taut. 

\section{The reduction principle}

We now discuss a reduction principle that will be used
to prove that the remaining four polar actions on $\mathbf{Ca}P^2$ are taut.
Let $G$ be a compact Lie group acting isometrically
on a connected, complete Riemannian manifold $X$. 
It is not necessary to assume that $G$ is connected.
Denote by $H$ a fixed
principal isotropy subgroup of the $G$-action on $X$, let $X^H$ be the
totally geodesic submanifold
of $X$ that is left point-wise fixed by the action of
$H$, and let ${}_cX$ denote the closure in $X$ 
of the subset of regular points of $X^H$. Then it is known
that ${}_cX$ consists of those components of $X^H$ 
which contain regular points of $X$, and that any component
of ${}_cX$ intersects all the orbits. 
Let $N$ be the normalizer of $H$ in $G$.
Then the group $\bar N=N/H$ acts on ${}_cX$ with trivial principal
isotropy subgroup in such a way that 
the inclusion ${}_cX\to X$ induces a homeomorphism
between orbit spaces ${}_cX/\bar N\to X/G$.
The pair $(\bar N,{}_cX)$ is called the \emph{reduction}
of $(G,X)$. In case ${}_cX$ is not connected, the pair
$(\bar N,{}_cX)$ can be further reduced as follows. 
Let $\bar X$ be a connected component of ${}_cX$, and let
$\bar G$ denote the subgroup of $\bar N$ that 
leaves $\bar X$ invariant. Then there is a homeomorphism
between orbit spaces $\bar X/\bar G\to{}_cX/\bar N$. 
We refer the reader to~\cite{Straume1,SS,G-S} for the proof
and a discussion of these assertions. 

From now on, we will work with a fixed
connected component $\bar X$ of $X^H$. 
Let $p\in\bar X$. It is clear that $Gp\cap\bar X=\bar Gp$. 
If $M=Gp$, we will denote $M\cap\bar X$ by $\bar M$.  
If, in addition, $p$ is a $G$-regular point, then 
the normal space to the principal orbit $M$ at~$p$
is contained in $T_p\bar X$, because the slice representation at $p$
is trivial.

There is a natural action of $G$ on the space 
$\mathcal P(X)$ of $H^1$-paths $\gamma:[0,1]\to X$
given by $(g\cdot\gamma)(t)=g\gamma(t)$, where
$g\in G$ and $t\in[0,1]$, since the Sobolev class $H^1$
is invariant under diffeomorphisms.  
In particular, 
if $q\in \bar X$ and $M$ is a $G$-orbit in $X$,
this action restricts to an action of $H$ on 
the space $\mathcal P(X,M\times q)$,
and it is clear that 
$\mathcal P(X,M\times q)^H=\mathcal P(X^H,M^H\times q)=
\mathcal P(\bar X,\bar M\times q)$.

\begin{lem}\label{critical-set}
Let $M$ be a $G$-orbit and $q\in \bar X$ be a
regular point for $G$. Consider the energy functional
$E_q:\mathcal P(X,M\times q)\to\mathbf R$.
Then the critical set of $E_q$ 
coincides with the critical set of the
restriction of $E_q$ to
$\mathcal P(\bar X,\bar M\times q)$.
\end{lem}

\Pf Let $\gamma$ be a critical point of 
$E_q:\mathcal P(X,M\times q)\to\mathbf R$.
Then $\gamma$ is a geodesic of $X$  
that is perpendicular to $M$ at $p=\gamma(0)$
and such that $\gamma(1)=q$.  
Since $\gamma$ is perpendicular to $M$ at $p$, it is
also perpendicular to $Gq$ at $q$. 
Since $Gq$ is a principal orbit, 
$G_q=H$ fixes $\gamma$ point-wise. 
It follows that $\gamma$ is a geodesic
of $X^H$, and thus of $\bar X$, that is perpendicular to $\bar Gq$,
from which we deduce that $\gamma$ is perpendicular
to $\bar Gp=\bar M$. 
Thus $\gamma$ is a critical point of 
$E_q:\mathcal P(\bar X,\bar M\times q)\to\mathbf R$.

On the other hand, let $\gamma$ 
be a critical point of $E_q:\mathcal P(\bar X,\bar M\times q)\to\R$. 
Then $\gamma$ is a geodesic of $\bar X$ that is perpendicular
to $\bar M=\bar Gp$ at $p=\gamma(0)$ and such that $\gamma(1)=q$. 
It follows that $\gamma$ is perpendicular to $\bar Gq$. 
Since $Gq$ is a principal orbit, 
the normal spaces of $Gq$ in $T_qX$ and of $\bar Gq$ in $T_q\bar X$
coincide. Therefore $\gamma$ is perpendicular
to $Gq$, and hence, to $Gp$. 
Since $\bar X$ is totally geodesic in $X$,
it follows that $\gamma$ is a critical point of $E_q:
\mathcal P(X,M\times q)\to\mathbf R$. \EPf

\medskip

Let $\mathcal Q$ be any subset of $\mathcal P(X)$. 
For $c>0$, denote by $\mathcal Q^c$ the 
subset $\mathcal Q\cap E^{-1}([0,c])$,
where $E:\mathcal P(X)\to\R$ is the energy functional.  
Also, denote by $\beta(Z)$ the (possibly
infinite) sum of the $\Z_2$-Betti numbers
of a topological space $Z$ (say,
with respect to singular homology). 

\begin{lem}\label{floyd}
Let $M$ be an arbitrary $G$-orbit,
$q\in\bar X$ and $c>0$. 
Suppose there is a subgroup $L\subset H$ which is a 
finitely iterated $\mathbf Z_2$-extension of the identity and such that 
the connected component of the fixed point set $X^L$
containing $q$ coincides with $\bar X$. 
Then the numbers $\beta(\mathcal P(\bar X,\bar M\times q)^c)$ and
$\beta(\mathcal P(X,M\times q)^c)$ are finite and 
\[ \beta(\mathcal P(\bar X,\bar M\times q)^c)\leq
\beta(\mathcal P(X,M\times q)^c). \]
\end{lem}

\Pf We will work with the finite 
dimensional approximations of the spaces
$\mathcal P(\bar X,\bar M\times q)^c$ and
$\mathcal P(X,M\times q)^c$~\cite{Mi}. Let $\rho$
denote the injectivity radius of $X$;
note that the injectivity radius of $\bar X$ 
cannot be smaller than $\rho$. Fix a positive integer
$n$ such that $n>c/\rho^2$. Denote by $\Omega$ 
the space of continuous, piecewise smooth curves $\gamma:[0,1]\to X$
of energy less than or equal to $c$
such that $\gamma(0)\in M$, $\gamma(1)=q$ and
$\gamma|[\frac{j-1}n,\frac jn]$ is a geodesic
of length less than $\rho$ for $j=1,\ldots,n$.
It is known that $\Omega$ is an open finite dimensional
manifold, and 
there is a deformation retract 
of $\mathcal P(X,M\times q)^c$ onto $\Omega$. 
It is clear that $\Omega$ is $H$-invariant;
moreover, $\mathcal P(\bar X,\bar M\times q)^c$
contains $\Omega^H$, and 
the definition of the deformation retract
shows that it 
restricts to a deformation retract 
of $\mathcal P(\bar X,\bar M\times q)^c$ onto
$\Omega^H$. 

We can write $L=\langle h_1,\ldots,h_n\rangle\subset H$
such that $h_{i+1}$ normalizes $H_i$ and $h_{i+1}^2 \in H_i$
for $i:0,\ldots,n-1$, where $H_0$ is the identity group 
and $H_i$ is the subgroup of $H$ generated by $h_1,\ldots,h_i$. 
Now each transformation $h_{i+1}$ has order $2$
in the fixed point set 
$\Omega^{H_i}$, 
so by a theorem of 
Floyd~\cite{Floyd}, 
\[ \beta(\Omega^L)\leq
\beta(\Omega^{H_{n-1}})\leq
\ldots\leq\beta(\Omega^{H_0}). \]
But $\Omega^L=\Omega^H$ and 
$\Omega^{H_0}=\Omega$,
so this completes the proof. \EPf

\begin{prop}\label{red}
Suppose there is a subgroup $L\subset H$ which is a 
finitely iterated $\mathbf Z_2$-extension of the identity and such that 
a connected component of the fixed point set $X^L$
coincides with $\bar X$. Suppose also
that the reduced action of 
$\bar G$ on $\bar X$ is taut.
Then the action of $G$ on $X$ is taut.
\end{prop}

\Pf Let $M$ be an arbitrary $G$-orbit. 
We will prove that $M$ is taut. 
Suppose $q\in X$ is not a focal point of $M$,
and consider the energy functional 
$E_q:\mathcal P(X,M\times q)\to \R$;
this is a Morse function. 
We must show that $E_q$ is perfect, and in fact
it is enough to do that for $q$ in a dense subset of $X$, so
we may take $q$ to be a regular point.
Since $G$ acts on $X$ by isometries, we may 
moreover assume that $q\in \bar X$. 
To simplify the notation, denote 
$\mathcal P(X,M\times q)$ by $\mathcal M$
and $\mathcal P(\bar X,\bar M\times q)$ by $\bar{\mathcal M}$.
If the action of $\bar G$ on $\bar X$ is taut, 
then each $\bar G$-orbit in $\bar X$ is taut,
hence $\bar M$ is taut. Note that $q$ is not a 
focal point of $\bar M$, so the restriction 
$E_q|\bar{\mathcal M}$ 
is a perfect Morse function. Therefore
$\mu_k(E_q|(\bar{\mathcal M})^c)=
\beta_k((\bar{\mathcal M})^c)$
for all $k$ and all $c$.
Let $\mu$ denote $\sum_k\mu_k$. Then  
\[ \mu(E_q|(\bar{\mathcal M})^c)=\beta((\bar{\mathcal M})^c), \]
for all $c$. 
Now Lemmas~\ref{critical-set} and~\ref{floyd} and the 
Morse inequalities imply that
\[ \beta((\bar{\mathcal M})^c)\leq\beta(\mathcal M^c)
\leq\mu(E_q|\mathcal M^c)=\mu(E_q|(\bar{\mathcal M})^c).\]
It follows that $\mu(E_q|\mathcal M^c)=\beta(\mathcal M^c)$
for all $c$. This implies that 
$\mu_k(E_q|\mathcal M^c)=\beta_k(\mathcal M^c)$
for all $k$ and all $c$, and thus 
$E_q$ is perfect. \EPf 

\medskip

We finish the proof of Theorem~\ref{main} 
by using Proposition~\ref{red}
to prove that the four polar actions 
on $X=\mathbf{Ca}P^2$ listed in~(\ref{gps})
are taut. For each action, it is enough 
to show that the reduced action is taut,
and to exhibit the group $L$ 
with the required properties. Some facts about the involved
actions that are used below can be found in~\cite{PTh1}.
In particular, it is shown there that the sections for
each one of these actions are isometric to $\R P^2$.  

\subsection*{$G=\Spin8$}

Since this group is contained in $\Spin9$,
it admits a fixed point $p\in X$. The slice representation
at $p$ is the sum of two inequivalent
$8$-dimensional representations of $\Spin8$.
A principal isotropy subgroup $H$ 
of $G$ on $X$ can be taken to be a principal isotropy
subgroup of the slice
representation, so $H\cong\G$. 
The homogeneous space $\Spin8/\G$ is diffeomorphic
to $S^7\times S^7$, so it follows that the normalizer $N$
of $\G$ in $\Spin8$ is isomorphic to $\Z_2^2\cdot\G$.  
Since $\bar N=N/H$ is discrete, $\bar G$ is also discrete,
so $\bar X$ is two-dimensional and thus
coincides with a section $\Sigma\cong\R P^2$. 
Now the reduced action of $\bar G$ on $\bar X$ is taut,
because its orbits are points and points are taut
in a symmetric space~\cite{TTh2}. The connected component 
containing $p$ of
the fixed point set 
of a subgroup $L$ of $H$ is completely determined by the
fixed point set of its linearized action on $T_pX$. 
The action of $H$ on $T_p\Sigma$ is of course trivial,
and its action on the orthogonal complement
of $T_p\Sigma$ in $T_pX$ 
is equivalent to the sum of two copies of 
the $7$-dimensional representation of $\G$.
Hence the required $L$ can be constructed as follows.
Consider the standard embedding of $\SU3$ in $\G$.
The normalizer of $\SU3$ in $\G$ is isomorphic to $\Z_2\cdot\SU3$;
let $h$ be a generator of the $\Z_2$-factor. We take
$L$ to be generated by the diagonal elements of $\SU3$
with $\pm1$ entries and by $h$. 

\subsection*{$G=\SO2\cdot\Spin7$}

There is also a fixed point $p$ in this case. 
The slice representation at $p$ is $\R^2\otimes\R^8$,
where $\R^2$ is the standard representation
of $\SO2$ and $\R^8$ is the spin representation of
$\Spin7$. The principal isotropy 
subgroup of $G$ on $X$ is
$H=\Z_2\cdot\SU3$. 
The normalizer $N$ of 
$H$ in $G$ is isomorphic to $\SO2\times\SO2\cdot\SU3$.
Since $\bar N=N/H$ is two-dimensional, $\bar X$ is a 
four-dimensional connected
complete totally geodesic submanifold of $X$ that
contains a section $\Sigma\cong\R P^2$.
It follows that $\bar X$ must be isometric to $\C P^2$
(compare Theorem~1 in~\cite{Wo}),
and then the reduced action of $\bar G$ 
on $\bar X$ is equivalent to the action of a maximal
torus $T^2$ of $\SU3$ on $\C P^2=\SU3/\mathbf{S}(\U1\times\U2)$. 
This action is known to be polar, and hence it is 
taut by the argument in section~\ref{lifting}. 
The subgroup $L$ of $H$ can be taken to be 
generated by the diagonal elements of $\SU3$
with $\pm1$ entries.

\subsection*{$G=\SU2\cdot\SU4$}

There is also a fixed point $p$ in this case. 
The slice representation at $p$ is $\C^2\otimes\C^4$,
where $\C^2$ is the standard representation
of $\SU2$ and $\C^4$ is the standard representation of
$\SU4$. Therefore the principal isotropy 
subgroup $H$ of $G$ on $X$ is $T^1\cdot\SU2$. 
The normalizer of $H$ in $G$ is $T^3\cdot\SU2$.
It follows that $\bar X$ is isometric to $\C P^2$,
and the reduced action is taut as in the preceding case.
The subgroup $L$ can be taken to be generated
by the diagonal elements of $H$
with $\pm1$ entries.

\subsection*{$G=\SU3\cdot\SU3$}

There is an orbit of type $\C P^2$. The 
corresponding isotropy is $\U2\cdot\SU3$, and 
the slice representation is $\C^2\otimes\C^3$.
Therefore the principal isotropy 
subgroup $H$ of $G$ on $X$ is $T^2$, 
the normalizer of $H$ in $G$ is $T^4$,
$\bar X$ is isometric to $\C P^2$,
and the reduced action is as taut in the preceding case.
The subgroup $L$ can be taken to be generated
by the diagonal elements of $H$
with $\pm1$ entries.


\providecommand{\bysame}{\leavevmode\hbox to3em{\hrulefill}\thinspace}
\providecommand{\MR}{\relax\ifhmode\unskip\space\fi MR }
\providecommand{\MRhref}[2]{%
  \href{http://www.ams.org/mathscinet-getitem?mr=#1}{#2}
}
\providecommand{\href}[2]{#2}

\end{document}